\documentclass[aos,nameyear,dvips]{arximspdf}

\doi{10.1214/10-AOS866}
\volume{38}
\issue{6}
\pubyear{2010}
\firstpage{3840}
\lastpage{3841}



\begin{document}
\begin{frontmatter}

\title{Letter to the Editor: Some comments on:\\
On construction of the smallest one-sided confidence interval
for the difference of two proportions}\vspace*{6pt}
\runtitle{Letter to the Editor}

\textit{Ann. Statist.} \textbf{38} (2010) 1227--1243

\begin{aug}
\author[A]{\fnms{Chris} \snm{Lloyd}\corref{}\ead[label=e1]{c.lloyd@mbs.edu}} and
\author[B]{\fnms{Paul} \snm{Kabaila}\ead[label=e2]{P.Kabaila@latrobe.edu.au}}
\runauthor{C. Lloyd and P. Kabaila}
\affiliation{Melbourne Business School and La Trobe University}
\address[A]{Melbourne Business School\\
Carlton, 3053\\
Australia \\
\printead{e1}}
\address[B]{School of Engineering \\
\hspace{1em} and Mathematical Sciences\\
La Trobe University\\
Victoria 3086\\
Australia\\
\printead{e2}}
\end{aug}

\received{\smonth{10} \syear{2010}}
\revised{\smonth{11} \syear{2010}}

%

%

\end{frontmatter}

\section{Introduction}

Wang (\citeyear{W10}) has rediscovered a general technique due to Buehler
[(\citeyear{B57}), Section 7] which was justified by Jobe and David
[(\citeyear{JD92}), Appendix A1] and then in more generality by \citet{LK03}.
The Buehler $1-\alpha$ upper confidence limit for a scalar parameter of
interest, based on a designated
statistic $L$, is $u(L)$ where $u$ is that nondecreasing function which
makes $u(L)$ as small as possible
subject to the constraint that the infimal coverage is $1-\alpha$.
Because Buehler illustrated the application of his result to the
reliability of a parallel system, his work
was virtually unknown outside the reliability literature for over 40
years. We believe that Buehler
confidence limits have many important statistical and computational
properties. The purpose of this letter is to point
the reader to some of the literature on these properties.

\section{Important statistical and computational properties of Buehler
confidence limits}

In the reliability literature on Buehler confidence limits, partly for
computational reasons,
the ordering induced on the sample space was usually based on an
estimator $L$ of the parameter of interest.
It turns out that this ordering typically leads to confidence limits
that do not have large sample efficiency;
see \citet{K01} and \citet{KL03}. To obtain large sample
efficiency, a Buehler $1-\alpha$
confidence limit needs to be based on an ordering induced on the sample
space by $L$ an approximate
$1-\alpha$ confidence limit. However, as noted by \citet{KL04}, such Buehler confidence
limits do not satisfy the nesting property. In the same paper, we
suggested a method of resolving the tension between
large sample efficiency and the satisfaction of the nesting property.
Of course, one seeks to obtain not only good
large sample performance, but also good finite sample performance.
Kabaila and Lloyd (\citeyear{KL02}, \citeyear{KL05}, \citeyear{KL06}) examine some of the
factors that influence the finite sample performance
of Buehler confidence limits. Proposition 2 of \citet{W10} is a
rediscovery of one result
in the latter paper. Buehler confidence limits also have computational
advantages which are examined
by \citet{K05}.

To summarize, Buehler confidence limits are exact, relatively easily
computed and possess an
attractive finite-sample optimality property. These advantages usually
come at the cost of either some
loss of large sample efficiency or nonsatisfaction of the nesting
property. Nonetheless,
Buehler confidence limits play an important part in statistical practice;
see, for instance, Lloyd and Moldovan (\citeyear{LM07a}, \citeyear{LM07b}).


%
\printaddresses

\end{document}